\documentclass{amsart}
\usepackage{amsmath}
\usepackage{amssymb}
\usepackage{amsthm}
\usepackage{amsfonts}

\newtheorem{lemma}{Lemma}
\newtheorem{definition}{Definition}

\newtheorem{theorem}{Theorem}
\newtheorem{proposition}{Proposition}

\theoremstyle{plain}
\newtheorem{thm}{Theorem}[section]

\theoremstyle{definition}

\theoremstyle{remark}
\newtheorem{rem}[thm]{Remark}

\def\Ua{U_{\alpha}}

\def\h{\textup{ht}}
\def\Id{\textup{Id}}
\def\not{/\!\!\!\!\!}
\def\rk{\textup{rk}}
\newcommand{\ra}{\to}
\newcommand{\Ph}[1]{(\Phi^{#1})'}
\newcommand{\ph}[1]{\Phi^{#1}}
\newcommand{\al}[1]{\alpha_{#1}}

\def\Id{\textup{Id }}

\newcommand{\epsa}[1]{\varepsilon_{\alpha_{#1}}}
\newcommand{\epsb}[1]{\varepsilon_{\beta_{#1}}}
\newcommand{\ignore}[1]{}
\newcommand{\field}[1]{\mathbb{#1}}

\newcommand{\F}{\field{F}}
\newcommand{\Z}{\field{Z}}
\newcommand{\beq}{\begin{displaymath}}
\newcommand{\eeq}{\end{displaymath}}
\newcommand{\set}[1]{\{#1\}}
\newcommand{\DD} {\mathcal{D}}

\begin{document}

\title[A diameter bound for unipotent groups]{A sharp diameter bound
 for unipotent groups of classical type over ${\mathbb{Z}}/p{\mathbb{Z}}$}
\author[Ellenberg and Tymoczko]{Jordan S. Ellenberg and Julianna
 Tymoczko}

\maketitle

\begin{abstract}  The unipotent subgroup of a finite group of Lie type
  over a prime field $\F_p$ comes equipped with a natural set of
  generators; the properties of the Cayley graph associated to this
  set of generators have been much studied.  In the present paper, we
  show that the diameter of this Cayley graph is bounded above and
  below by constant multiples of $np + n^2 \log p$, where $n$ is the
  rank of the associated Lie group.  This generalizes the result of
  \cite{elle:thesis}, which treated the case of $SL_n(\F_p)$.
  (Keywords:  diameter, Cayley graph, finite groups of Lie type.  AMS
  classification:  20G40, 05C25)
\end{abstract}

\section{Introduction}\label{intro}

Given a group $G$ endowed with a set of generators $\Sigma = \{u_1,
\ldots, u_k\}$, the {\em Cayley graph} of $G$ in $\Sigma$ is the
directed graph whose vertices are elements of $\Sigma$, and in which
$x$ and $y$ are joined if $y = u_i x$ for some $i$.  It is a
long-standing problem to investigate geometric properties of the
Cayley graph; especially interesting is the behavior of the random
walk obtained by setting off along a randomly chosen edge at each
step.  In case $G$ is a finite group, a natural invariant of the
Cayley graph of $G$ in $\Sigma$ is its {\em diameter} -- that is, the
minimum positive integer $N$ such that every element of $G$ can be
written as a word of length at most $N$ in $\Sigma$.  If $g$ is an
element of $G$, we define the {\em length} of $g$ in $\Sigma$ (or just
the length of $g$ if $\Sigma$ is understood) to be the length of the
shortest word in $\Sigma$ which evaluates to $g$.

Diameters of Cayley graphs have been much studied: for a general
discussion, see \cite{bhkls} or \cite[\S3.8]{b1}.  An excellent
reference for the theory of random walks on finite groups in general
is Saloff-Coste's survey paper \cite{salo:rwfg}.  The relationship
between the diameter of $G$ in $\Sigma$ and the convergence of the random
walk along $\Sigma$ is discussed, for example, in \cite[Corollary
1]{dsc93}; for random walks on unipotent groups (in particular, groups
of upper-triangular matrices) see \cite{cp}, \cite{ads}, and
\cite{stong}.  The problem of controlling the diameters in natural
families of groups is quite hard.  For a general permutation group
with a general set of generators, the problem of computing the
diameter is NP-hard \cite{eg}.  Bounding the diameter of $SL_2(\F_p)$
in standard generators, as $p$ grows, is a well-known problem whose
solution (at present) requires the theory of automorphic forms (but
see \cite{helf} for a new approach via additive number theory.)

In this paper, we give upper and lower bounds for the diameter of
unipotent subgroups of classical groups over prime fields; such
finite groups admit a natural choice of $\Sigma$, which we describe in
section~\ref{Chev}.  In case the classical group is $SL_n$, we recover
the previously unpublished result of the first
author~\cite{elle:thesis} cited in \cite{ads}, \cite{dsc94}, and
\cite{salo:rwfg}.  The main theorem of our paper is the following:

\begin{thm}
There exists absolute constants $c$ and $C$ such that, if $p$
  is an odd prime, and $G$ is a Chevalley group of classical type of rank $n$ over $\Z/p\Z$, the diameter of $U$ in $\Sigma$ is at most $C(np + n^2 \log
  p)$ and at least $c(np + n^2 \log p)$.
\label{main theorem}
\end{thm}

In the first section of this paper, we outline some of the basic facts
we will need about the combinatorics of commutator relations in
unipotent groups of classical type; in the second and third sections,
respectively, we prove the upper and lower bounds on the diameter of
$U$ in $\Sigma$.  In an appendix we provide more detailed descriptions
of the Chevalley groups as matrix groups.
ecomposing along the ``rows'' of the matrix groups.  

To begin, we summarize the algebra needed for this task.  
Subsection \ref{Chev}
gives a brief introduction to Chevalley groups and Subsection \ref{rows} 
details a decomposition of $G$ upon which we rely heavily.  Subsection
\ref{generators} uses both to give generators for the subgroups we will
study. 

{\bf Acknowledgments.}  The authors are grateful to Persi Diaconis,
under whose supervision the first author completed a special case of
the present work; and to Igor Pak, for many helpful comments and
suggestions.  The first author is partially supported by NSF Grant
DMS-0401616.  The second author is partially supported by an 
NSF Postdoctoral Fellowship and received support from 
the Clay Mathematics Institute during the writing of this paper.

\subsection{Properties of Chevalley groups of classical
  type}\label{Chev}
Let $K$ be a field of characteristic other than $2$, and let $G$ be a
Chevalley group of classical type over the field $K$: that is, $G$ is
either a special linear group, a special orthogonal group, or a
symplectic group over $K$.  Because $G$ is a Chevalley
group, its multiplicative structure is described combinatorially by a
set of roots and the relations between these roots.  In this section
we review the combinatorics associated to $G$ and some basic results
we will need for our arguments.  The Appendix provides explicit matrix
representatives for all these constructions.

Let $B$ be a Borel subgroup of $G$, and let $U$ be the unipotent
subgroup of $B$; so $B$ factors as a product $TU$ where $T$, the
maximal torus contained in $B$, is isomorphic to $(K^*)^n$ for an
integer $n$ called the rank of $G$.  For more details about the
structure of Chevalley groups, see \cite[Chapters 8, 28, and 35]{H1}.

In this paper, we will give upper and lower bounds for the diameter of
$U$ in a natural set of generators described below, in the case where
$K$ is $\Z/p\Z$ for some odd prime $p$. 

We begin by reviewing several properties of $U$.  First, $U$ is
generated by a collection of one-parameter subgroups $U_{\alpha}
=\{\epsa{}(c): c \in K\}$ called {\em root subgroups}; here $\alpha$
ranges over the roots of $G$, and the group homomorphism $\epsa{}: K
\longrightarrow G$ is obtained by exponentiating the eigenspace of the
Lie algebra of $G$ corresponding to $\alpha$.  The set of $\alpha$
such that $\epsa{}(K) \subset U$ is called the set of positive roots,
and is denoted $\Phi^+$.


The positive roots 
inherit an operation from the character group 
which we denote addition.  (Warning: the set of positive roots is not
closed under this operation!)  A minimal generating set of $\Phi^+$
will have cardinality $n$; denote such a generating set by
$\Delta = \{\alpha_i: 1 \leq i \leq n\}$.  These $\alpha_i$ are
called the simple roots.  For example, in type $A_n$ (i.e. when $G =
SL_{n+1}$) the subgroup $U$ is the group of upper-triangular matrices with
$1$ along the diagonal while the simple roots correspond to the root
subgroups $U_{\alpha_i} = \{\Id + cE_{i,i+1}: c \in K\}$.  More
examples of explicit matrix representatives are in the Appendix.  We
use the conventions of \cite{H2} to index the simple roots; Figure
\ref{roots} describes all the positive roots of $G$ with respect to
this choice of indexing.

The set $\Sigma \subset G$ of elements of the form $\epsa{i}(\pm 1)$
is a generating set for $G$ of cardinality $2n$.  The first part of
this paper will be devoted to proving the following proposition:

\begin{proposition} There exists an absolute constant $C$ such that, if $p$
  is an odd prime, and $G$ is a Chevalley group of classical type over
  $Z/pZ$, the diameter of $U$ in $\Sigma$ is at most $C(np + n^2 \log
  p)$.
\label{upper bound}
\end{proposition}

\begin{rem} The constant $C$ is easy to bound explicitly, if desired.
\end{rem}

\begin{rem} The difficult part of the argument is controlling the
  behavior of the diameter of $G$ as $n$ grows; this is why we do not
  consider diameters of exceptional groups of Lie type, in which only
  the growth of $p$ is at issue.
\end{rem}

The key idea driving our argument is that the non-commutativity of
$U$ can be described in completely combinatorial terms by means of the
Chevalley commutator relation, which is the following:

\begin{equation}
[ \epsa{}(s), \epsb{}(t) ] = \prod_{\textup{$ 
    \begin{array}{c}i,j > 0 \vspace{-.5em} \\ 
       i \alpha + j \beta \in \Phi^+\end{array}$}} 
        \varepsilon_{i \alpha + j \beta } (c_{ij} s^i t^j)
\label{ccr}
\end{equation}
(see \cite[page 207]{S}).
Brackets denote the group commutator $[x,y]=xyx^{-1}y^{-1}$ in $U$.  
The integers $c_{ij}$ are 
determined by the roots $\alpha$ and $\beta$ as well as the ordering of 
the product; we list many of them in the Appendix.  The following proposition
is the only fact we will use about the $c_{ij}$.

\begin{proposition}
If $\textup{char } K 
\neq 2$ then $c_{ij} \neq 0$ if and only if $i \alpha + j \beta
\in \Phi^+$.
\end{proposition}

The additive relations between roots also play a key role in 
our computations. Figure \ref{roots} enumerates
the positive roots of classical Chevalley groups, for $i$ between $1$ and
$n = rk(G)$.

\begin{figure}[h]
 \begin{tabular}{|c|c|c|}
        \cline{1-3} Root & Parameters & Type \\
        \cline{1-3}       
        $ \sum_{j=i}^k \alpha_{j} $ &
        $1 \leq i \leq k \leq n$ &
        $A_n, B_n, C_n, D_n$ \\
        \cline{1-3}\multicolumn{3}{|c|}{\textup{except 
          $\alpha_{n-1}+\alpha_n \notin \Phi$ in type $D_n$}} \\
        \cline{1-3}
        $ \sum_{j=i}^n \alpha_{j} + \sum_{j=k}^n \alpha_j$ &
        $1 \leq i < k \leq n$ &
        $B_n$ \\
        \cline{1-3}
        $ \sum_{j=i}^n \alpha_{j} + \sum_{j=k}^{n-1} \alpha_j$ &
        $1 \leq i \leq k < n$ &
        $C_n$ \\
        \cline{1-3}
        $\sum_{j=i}^{n-2}\alpha_{j} + \alpha_n$ &
        $1 \leq i \leq n-2$ &
        $D_n$ \\        
        \cline{1-3}
        $ \sum_{j=i}^{n} \alpha_{j}  + 
        \sum_{j=k}^{n-2} \alpha_{j}$ &
        $1 \leq i < k \leq n-2$ &
        $D_n$ \\
        \hline
\end{tabular}
\caption{Positive roots of groups of classical type} \label{roots}
\end{figure}
%

As we can see in Figure \ref{roots}, each positive root can be written
uniquely as a linear combination of simple roots with nonnegative
integer coefficients.  The sum of these coefficients is called the
height of the root and is denoted $\h(\sum c_i \alpha_i) = \sum c_i$.
We write $\alpha > \beta$ if $\alpha - \beta$ is a sum of positive
roots. 

For a thorough introduction to algebraic groups and root systems, the 
reader is referred to \cite{S} or \cite{H1}.

\subsection{Decomposing $U$ into $U_i$} \label{rows}
To describe the diameter of $U$ we decompose $U$ as a product of
subgroups $U_i$, by analogy with the decomposition of the group of
upper triangular matrices as a product of the subgroups consisting of
matrices with zero entries away from the diagonal and row $i$.  In
this section we define the $U_i$ and prove some of their properties.
In later sections we will bound the diameter of
each $U_i$ using the Chevalley commutator relation.

We begin by defining a subset $\Ph{i}$ of $\Phi^+$:
\[\Ph{i} =  \{ \al{i}\} \cup
        \{ \al{} \in \Phi^+ | \al{} > \alpha_i, /\!\!\! \exists j < i \textup{ such that } \alpha > \alpha_j\}. \]
For instance, in type $A_3$ these sets are $\Ph{1} = \{\alpha_1, \alpha_1 + 
\alpha_2, \alpha_1 + \alpha_2 + \alpha_3 \}$, $\Ph{2} = \{ 
\alpha_2, \alpha_2 + \alpha_3 \}$, and $\Ph{3} = \{\alpha_3 \}$.  By contrast,
in type 
$D_5$ the set $\Ph{2}$ is $\Ph{2} = \{\alpha_2, \alpha_2 + \alpha_3, 
\alpha_2 + \alpha_3 + \alpha_4, \alpha_2 + \alpha_3 + \alpha_5, 
\alpha_2 + \alpha_3 + \alpha_4 + \alpha_5, \alpha_2 + 2 \alpha_3 + \alpha_4 + 
\alpha_5\}$.  Note that, for fixed $i$, the roots given in the first
column of Figure \ref{roots} are precisely those contained in $\Ph{i}$.

Let $U_i'$ be the subgroup of $U$ generated by 
$\{U_{\alpha}: \alpha \in \Ph{i}\}$.  The Chevalley relations 
and Figure \ref{roots} show that $U_i'$ is 
abelian in types $A_n$,$B_n$,
and $D_n$, and Heisenberg in type $C_n$.  In other words, since the
set $\{\alpha+\beta: \alpha \in \Ph{i}, \beta \in \Ph{i}, 
\alpha+\beta \in \Phi^+\}$ is empty in types 
$A_n$, $B_n$, and 
$D_n$, the elements of $U_i'$
commute.  However, in type $C_n$ 
there is a unique root $\gamma_i = 2 \sum_{j=i}^{n-1} 
\alpha_j + \alpha_n$ that has the property
\[\{\gamma_i\} = 
\{\alpha+\beta: \alpha \in \Ph{i}, \beta \in \Ph{i}, 
\alpha+\beta \in \Phi^+\}.\]  
Moreover, for 
each $\alpha \neq \gamma_i$ in $\Ph{i}$ the root
 $\gamma_i - \alpha$ is
also in $\Ph{i}$.  The root $\gamma_i$ is called the long root in $\Ph{i}$.

The roots in $\Ph{i}$ in type $A_n$ have the form 
$\sum_{j=i}^k \alpha_j$ for $k = i, i+1, \ldots, n$ and so
are totally ordered by height.  In fact, if the heights of 
$\alpha$ and $\beta$ in $\Ph{i}$ differ by one then 
$\alpha - \beta= \pm \alpha_j$ for some simple root $\alpha_j$.
Moreover, if $\alpha$ is in $\Ph{i}$ and $\alpha + \alpha_j$ is
a root then $\alpha + k \alpha_j$ is not a root for any $k>1$.

For a general Chevalley group $G$, we will define $\ph{i}$ and $U_i$ so as
to preserve these properties to the greatest extent possible.  Our
definitions are given in Figure \ref{phi}.  For instance, in type
$A_n$, the group $U_i$ consists of matrices $\{(\Id + \sum_{j=i+1}^n
a_j E_{ij}): a_j \in K\}$.  The Appendix gives detailed descriptions
of the $U_i'$ in the other types.

Note that in types $B_n$, $C_n$, and $D_n$, the $U_i$ are abelian
quotient groups of $U'_i$.   The images of the root subgroups $\Ua$
for $\{\alpha \in \Phi^i\}$ in the quotient $U_i' \longrightarrow U_i$
generate $U_i$.  By an abuse of notation, we also use
$\Ua$ and $\epsa{}$ to denote the image of $\Ua$ and $\epsa{}$ under
the homomorphism $U'_i \ra U_i$.

\begin{figure}[h]
\begin{tabular}{c|c|c} 
Type & $U_i$ & $\Phi^i$ \\
\cline{1-3} $A_n$ & $U_i'$ & $\Ph{i}$  \\
$B_n$ & $U_i'/U_{\sum_{j=i}^n \alpha_j}$ & $\Ph{i} - 
         \{\sum_{j=i}^n \alpha_j\}$ \\ 
$C_n$ & $U_i'/U_{\gamma_i}$ & $\Ph{i} - \{\gamma_i\}$ \\
$D_n$ & $U_i'/U_{\sum_{j=i}^{n-2} \alpha_j + \alpha_n}$ &
           $\Ph{i} - \{\sum_{j=i}^{n-2} \alpha_j + \alpha_n\}$
\end{tabular}
\caption{Definition of $U_i$ and $\ph{i}$} \label{phi}
\end{figure}

For all Chevalley groups of classical type, the group $U$ is the product
\[ U = \prod_{i=1}^n U_i' \]
\cite[11.3.3 Exercise 2a]{S}.
If $g = \prod_{\alpha \in \Phi^i} \epsa{}(s_{\alpha})$ and $g'= \prod_{\alpha
\in \Phi^i} \epsa{}(s_{\alpha}')$ are in $U_i$ then by the commutativity
of $U_i$ their product is
\begin{equation}\label{product in row}
gg' = \prod_{\alpha \in \Phi^i} \epsa{}(s_{\alpha} + s_{\alpha}').
\end{equation}

The next fact motivates the choice we have made of $\ph{i}$.
It follows by inspection of Figure \ref{roots}.  For each type, define
$r_i$ to be the cardinality of $\Phi^i$.

\begin{proposition} \label{total order}
%
%
$\ph{i}$ is totally ordered by height.  If $\alpha$ is in $\ph{i}$ 
and $\beta$ is a simple root 
then $j \alpha + l \beta$ is in $\ph{i}$ for exactly one
positive pair $(j,l)$, with $j=1$ and $l \in
\{ 1,2\}$.

In types $A_n$, $C_n$, and $D_n$
the set $\ph{i}$ contains exactly
one root of each height from $1$ to $r_i$.  In type $B_n$ the set
$\ph{i}$ contains exactly one root of each height from $1$ to $n-i$ and
one root of each height from $n-i+2$ to $r_i+1$.  
\end{proposition}

\subsection{Generators for $U_i$} \label{generators}

In the rest of this section, we assume that $i$ has been fixed, that 
$\ph{i}$ has been ordered by height from smallest to largest, and that
$\beta_j$ denotes the $j^{th}$ element of $\ph{i}$ in this order.

For instance, in type $A_n$ the matrices $\epsb{j}(t_j)$ have $t_j$ in
the $(i,j+i)$ entry, $1$'s along the diagonal, and zeroes elsewhere.

The following proposition uses the Chevalley commutator relation to
generate elements of large-height root subgroups by means of
reasonably short words. 

\begin{proposition} \label{commutator generation}
Define functions 
\[\begin{array}{l}
f:\{1,2,\ldots,r_i\} \longrightarrow
\{i, i+1,\ldots,\rk(G)\} \textup{ and} \\
m: \{1,2,\ldots, r_i\} \longrightarrow \{1,2\} \end{array}\]
by the equations 
$\beta_j-\beta_{j-1} = m(j) \alpha_{f(j)}$ (when $j \geq 2$), as well
as $m(1)=1$ and $f(1)=i$.

If $\beta_{j-1} \in \ph{i}$ and $l \in \{1, 2, \ldots, n\}$ then
\[ [\epsb{j-1}(s), \epsa{l}(t)] = \left\{ \begin{array}{ll}
     \epsb{j}(c_{1,m(j)} s t^{m(j)}) & \textup{ if } f(j)=l \\
     \epsb{j}(0) & \textup{ if } f(j) \neq l.  \end{array} \right.\]
\end{proposition}

Here $c_{1,m(j)}$ is the structure constant arising in Equation
\eqref{ccr}.
Note that the equality in Proposition~\ref{commutator generation}
holds in $U_i$; the corresponding equality in $U_i'$ does not
hold in general.  However, the equality {\em does} hold in $U_i'$ when $j
\leq n-i$.

\begin{proof}
The functions $f$ and $m$ are well defined 
by inspection of the tables in Figures \ref{roots} and \ref{phi}.  
In fact, note
that $m(j)=2$ only when $j=n-i+1$ in type $B_n$.  

By Proposition \ref{total order}, 
there is at most one root in $\ph{i}$ that is a 
linear combination of $\beta_{j-1}$ and $\alpha_l$ with positive integer 
coefficients.  If it exists, this root is $\beta_j$.

The claim is now an application of the Chevalley commutator relation.
\end{proof}

This proposition permits us to define a map from $K^r$ to $U_i$ that
generates elements of $U_i$ via successive conjugation by an element of
a simple root subgroup.  It is this map that allows us quickly to
generate elements of $U$ with large matrix entries.

\begin{definition}
Define the map $\theta_{k,r}: K^r \longrightarrow U_i$
inductively by
\[\begin{array}{l}
  \theta_{k,1}(s) = \epsb{k}(s) \\
  \theta_{k,r}(s_1, \ldots, s_r) = \left(\epsa{f(k+r-1)}(-s_r) \right)
  \left(\theta_{k,r-1}(s_1, \ldots, s_{r-1})\right) 
   \left(\epsa{f(k+r-1)}(s_r) \right) \end{array}\]
\end{definition}

$\theta_{k,r}$ is well-defined only if $k+r-1 \leq r_i$.
For example, if $n \geq 3$ in any classical type then the map $\theta_{1,3}$
is
\[
\theta_{1,3}(s_1, s_2, s_3) = \epsa{3}(-s_3) 
   \left( \epsa{2}(-s_2) \left( \epsa{1}(s_1) \right) 
    \epsa{2}(s_2) \right) \epsa{3}(s_3).\]
The matrix for this product when $i = 1$ in type $A_3$ (i.e. when $G = SL_4$) is
$\left(\begin{array}{cccc} 1 & s_1 & s_1 s_2 & s_1 s_2 s_3 \\
   0 & 1 & 0 & 0 \\ 0 & 0 & 1 & 0 \\ 0 & 0 & 0 & 1 \end{array}\right)$.

Note that $\theta_{k,r}(s_1, \ldots, s_r) \theta_{k,r}(-s_1, s_2, s_3,
\ldots, s_r)  = 1$.

Since the image of $\theta_{k,r}$ is in $U_i$ it factors uniquely as a 
product of elements of root subgroups, described next.

\begin{proposition}\label{thetaprops}
Write $\theta_{k,r}(s_1,
\ldots, s_r) = \prod_{l=1}^{r_i} \epsb{l}(t_{l})$.  
\begin{enumerate}
\item If $l < k$ or $k+r-1 < l$ then $t_l = 0$. \label{zeros}
\item Fix $k$ and $r$ so that either $n-i+1 \leq k$ or $k+r-1 \leq n-i+1$.
If $k \leq l \leq k+r-1$ 
then $t_l = s_k \prod_{j=k+1}^{l} c_{1,m(j)} s_{j}^{m(j)}$. \label{factors}
\end{enumerate}
\end{proposition}

\begin{proof}
The proof inducts on $r$.  The claim holds by definition if $r=1$.  Note
\[
\theta_{k,r}(s_1, \ldots, s_r) = \left(\epsa{k+r-1}(-s_{r}) \right)
   \left( \theta_{k,r-1}(s_1, \ldots, s_{r-1}) \right) 
    \left(\epsa{k+r-1}(s_{r}) \right)\]
by construction.  Use induction to write
$\theta_{k,r-1}(s_1, \ldots, s_{r-1}) = \prod_{l=k}^{k+r-2}
\epsb{l}(t_l')$.  Proposition \ref{commutator generation} shows that
\[\begin{array}{rl}
\displaystyle \left(\epsa{k+r-1}(-s_{r}) \right) \hspace{-1em}&
 \displaystyle  \left( \prod_{l=k}^{k+r-2}
      \epsb{l}(t_l') \right) \left(\epsa{k+r-1}(s_{r}) \right) 
   \\ &= \displaystyle \left( \prod_{l=k}^{k+r-2}
      \epsb{l}(t_l') \right) \prod_{{\scriptsize{\begin{array}{c}
            l \textup{ s.t. } k \leq l-1 \leq k+r-2 \\
         \textup{and } f(l) = k+r-1    \end{array}}}} 
        \epsb{l}(
        c_{1,m(l)} t_{l-1}' s_{r}^{m(l)}).\end{array}\]
This can be written as $\prod_{l=k}^{k+r-1} \epsb{l}(t_l)$ for
some $t_l$.  In particular, Part \ref{zeros} holds.

By inspection of Figure \ref{roots}, there are 
at most two $\beta_l$ in $\ph{i}$ with
$\beta_l - \beta_{l-1} = m(l) \alpha_{k+r-1}$, of which 
one has $l-1 < n-i+1$ while the other has $l-1 \geq n-i+1$.  By
hypothesis in Part \ref{factors}, at most one of these roots is in
$\{\beta_k, \ldots, \beta_{k+r-1}\}$.  
The root $\beta_{k+r-1}$ is such a root by 
definition of $\theta_{k,r}$.  So $t_l = t_l'$ for $l=k, \ldots, k+r-2$ 
and $t_{k+r-1} = c_{1,m(k+r-1)}
t_{k+r-2} s_r^{m(k+r-1)}$.  Induction completes the proof.
\end{proof}

\section{Proof of Proposition \ref{upper bound}}

In this section, we will prove Proposition \ref{upper bound}.  The
main idea is to show explicitly how to generate elements of $U_i$ with
relatively short words in $\Sigma$.  We accomplish this by further
decomposing each $U_i$ into three smaller subgroups, each generated by
root subgroups whose heights lie in a specified range.

\subsection{The contribution of the simple root subgroups}

\label{abel}

Each $U_i$ contains the simple root subgroup $U_{\alpha_i}$.  We can
generate an element $\epsa{i}(k)$ of $U_{\alpha_i}$ by means of at
most $p/2$ copies of $\epsa{i}(\pm 1)$.

\subsection{The contribution of medium-height root subgroups} \label{first k}

Fix $k \leq r_i$.  Let $U_i^k$ be the subgroup $\prod_{j=2}^k 
U_{\beta_j}$ of $U_i$ and let $U_i^{>k}$ be the subgroup $U_i^{>k} =
\prod_{j=k+1}^{r_i} U_{\beta_j}$. We study $U_i^k$ in this section and
$U_i^{>k}$  in the next.

The following lemma expresses each
generator of the root subgroup $U_{\beta_l}$ in terms of $\Sigma$.
Recall that $c_{1,m(j)}$ is the structure constant defined by
$[\epsb{j-1}(1), \epsa{f(j)}(1)]$ in $U_i$, as in
Proposition~\ref{commutator generation}.  
Denote the product of the first $l$ structure constants by
\[d_l = \prod_{j=2}^l c_{1,m(j)}.\]

(Note that $|d_l|$ is a power of $2$.)

\begin{lemma}\label{maxlengthmiddlegenerators}
For each $\beta_l$ in $\Phi^{i}$, the generator $\epsb{l}\left( d_l \right)$
of $U_{\beta_l}$ in $U_i$ can be written as a word in $\epsa{j}(\pm
1)$ with at most $8l$ letters. 
\end{lemma}

\begin{proof}
We write ${\bf 1^l}$ for the vector in $K^l$ each of whose entries is $1$.

First, assume $l \leq n-i+1$ 
or, equivalently,  $\beta_l \leq \sum_{j=i}^{n-1} \alpha_j + 2\alpha_n$.  
(The sum $\sum_{j=i}^{n-1} \alpha_j + 2\alpha_n$ is a root only in type $B_n$,
where it equals $\beta_{n-i+1}$.  In the other types, the condition
$\beta_l \leq \sum_{j=i}^{n-1} \alpha_j + 2\alpha_n$ is equivalent to 
$\beta_l \leq \sum_{j=i}^n \alpha_j$.)
From Proposition \ref{thetaprops} we have that
\begin{eqnarray*}
\theta_{1,l}({\bf 1^l}) 
  \theta_{1,l-1}(-1, {\bf 1^{l-2}}) &=&
  \theta_{1,l-1}({\bf 1^{l-1}}) \theta_{1,l-1}(-1, {\bf 1^{l-2}}) 
   \epsb{l}\left( d_l \right) \\
 &=& \epsb{l}\left( d_l \right).\end{eqnarray*}
By construction, the left hand side is a word of length at most $4l$ 
in the letters $\epsa{j}(\pm 1)$.  The element $\epsb{l}(d_l)$ gives a generator for $U_{\beta_l}$ in $U_i$ since $p$ is odd.  Note that $\epsb{n-i+1}\left(- d_{n-i+1} \right) = (\epsb{n-i+1}\left( d_{n-i+1} \right))^{-1}$ also has length at most $4(n-i+1)$.


Now suppose $l > n-i+1$.  The generator $\epsb{l}\left( d_l \right)$ 
can be written as
\[\theta_{n-i+1,l-n+i}\left(d_{n-i+1} ,
{\bf 1^{l-n+i-1}}\right) 
\theta_{n-i+1,l-n+i-1}\left(-d_{n-i+1},
{\bf 1^{l-n+i-2}}\right).\]
The length of this word is at most $8(n-i+1) + 4(l-n+i-1) \leq 8l$. 
\end{proof}

We now bound the length of $\epsb{l}\left(s d_l \right)$ 
for arbitrary $s$ in $K$.  

\begin{lemma} \label{rootelements}
The length of $\epsb{l}\left(s d_l \right)$ in $U_i$ 
is at most $48l\sqrt{s}$.
\end{lemma}

\begin{proof}
By abuse of notation we also use $s$ to denote the integer in the
real interval $[0,p-1]$ which reduces to $s$ mod $p$.  Let $t$ be the largest
integer less than or equal to $\sqrt{s}$:  
\[ t = \lfloor \sqrt{s} \rfloor. \]
Since the $\epsb{l}$ are additive homomorphisms, we know that
\[\epsb{l}\left(sd_l \right) 
   = \epsb{l}\left(t^2 d_l \right)
   \epsb{l}\left((s-t^2)d_l \right).\] 
By the commutator relationship in $U_i$, 
\[\begin{array}{lcl}
\vspace{.25em} \epsb{l}\left(t^2 d_l \right) 
&=& \left[\epsb{l-1}\left(t d_{l-1} \right),\epsa{f(l)}(t) \right] 
    \textup{ and}\\
\epsb{l}\left((s-t^2) d_l \right) &=& 
\left[\epsb{l-1}\left( d_{l-1} \right),\epsa{f(l)}(s-t^2)
\right].
\end{array}\]
When $l=n-i+1$ in type $B_n$, we use the relations
\[\begin{array}{lcl}
\vspace{.25em} \epsb{l}\left(t^2 d_l \right) &=& 
\left[\epsb{l-1}\left( d_{l-1} \right),
   \epsa{f(l)}(t)\right] \textup{ and}\\
\epsb{l}\left((s-t^2) d_l \right) &=& 
\left[\epsb{l-1}\left((s-t^2) d_l \right),
  \epsa{f(l)}(1)\right].\end{array}\]
Using the inequalities
\[s-t^2 \leq ((t+1)^2 - 1) - t^2 \leq 2t \leq 2\sqrt{s},\]
and the previous lemma 
we conclude that the elements 
$\epsb{l-1}\left(t d_{l-1} \right)$ and 
$\epsb{l-1}\left(d_{l-1} \right)$ 
have length at most 
$(8l-8)\sqrt{s}$ and $8l-8$.  In the type $B_n$ calculation, note
that the element 
$\epsb{l-1}\left((s-t^2) d_{l-1} \right)$ has length
at most $(8l-8)(2 \sqrt{s})$.  In all types, the elements
$\epsa{f(l)}(t)$ and $\epsa{f(l)}(s-t^2)$ have length at most $\sqrt{s}$ and
$2\sqrt{s}$, respectively.  

Consequently, the total length of 
$\epsb{l}\left(s d_l \right)$ is at most 
$48l \sqrt{s}$.
\end{proof}

\begin{proposition} \label{boundingUk}
Let $k = \left\lfloor \frac{\log p}{\log 2} \right\rfloor$.  
Every element of $U_i^k$ can be expressed by a word in $\Sigma$ of
length at most
\[\frac{48}{\log^2{2}}\sqrt{p} \log^2{p}.\]
\end{proposition}

\begin{proof}
Since $s$ is at most $p$, the diameter of $U_i^k$ is at most
\[\sum_{l=2}^k 48l\sqrt{p} < 48k^2\sqrt{p} \leq 
         \frac{48}{\log^2{2}}\sqrt{p} \log^2{p}.\]
\end{proof}

Note that the constant here is not meant to be optimal; in particular,
one can do better by specifying the type of the Chevalley group.


\subsection{The contribution of large-height root subgroups}\label{last k}

The next step is to bound the diameter of $U_i^{>k}$, which we do using
a kind of binary expansion as described below.

Fix an element $g$ in $U_i^{>k}$ and write 
$g =  \prod_{l=k+1}^{r_i} \epsb{l}(s_{l})$. 
Define the vector
\[
  {\bf u} = 
   (s_{k+1},s_{k+2},\ldots,s_{r_i}) \in K^{r_i-k}. 
\]

Recall that $d_l = \prod_{j=2}^l c_{1,m(j)}$ is the product of structure
constants.  
Each entry $\frac{s_l}{2 d_l}$ is in $K$ and so (considered as an
integer between $0$ and $p-1$) has a binary decomposition. 
Define a function $b(l,j)$ by
\[b(l,j) = \left\{ \begin{array}{ll} 0 &\mbox{  if $2^j$ is not in the
                binary decomposition of $\frac{s_l}{2 d_l}$} \\
                1 &\mbox{  if $2^j$ is in the
                binary decomposition of 
                $\frac{s_l}{2 d_l}$.} \end{array}
                \right. \]

Define
\[
  {\bf u_j} = 2^{j+1} \left(b(k+1,j) d_{k+1}, 
  \ldots, b(r_i,j) d_{r_i} \right). 
\]
so that $\sum_j {\bf u_j} = {\bf u}$.  We may think of the ${\bf u_j}$ as
making up a ``binary decomposition'' of {\bf u}.  If
$s_{j,l}$ denotes the $l^{th}$ entry of ${\bf u_j}$, then
\[\prod_{j=0}^{\left\lfloor \frac{\log p}{\log 2} \right\rfloor} \prod_{l=k+1}^{r_i} \epsb{l}(s_{j,l}) = 
\prod_{l=k+1}^{r_i} \epsb{l}(s_{l})=g.\]
Thus the length of $g$ is bounded by
the combined lengths of the words $\prod_{l=k+1}^{r_i} \epsb{l}(s_{j,l})$.
%


We first study words whose
vectors ${\bf u_j}$ are zero in the last $r_i-n+i-1$ entries.

\begin{lemma} 
Let ${\bf u} \in K^{n-i+1-j}$ be a vector with entries
$s_l \in \{0,d_{l+j}2^{j+1}\}$.
The word $\prod_{l=j+1}^{n-i+1} \epsb{l}(s_{l-j})$
has length at most $8(n-i+1)$.
\label{beforemiddle}
\end{lemma}

\begin{proof}
Define vectors ${\bf a_1}$ and ${\bf a_2}$ by
\[a_{1,l} = \left\{ \begin{array}{ll}
                        2 & 1 \leq l \leq j \\
                        1 &  j < l \leq n-i+1 \end{array} \right. \]
and
\[a_{2,l} = \left\{ \begin{array}{ll}
                        -2 & l=1 \\
                         2 &  1 < l \leq j \\
                        \left( -1 \right)^{\frac{s_{l-j}}{d_{l} 2^{j+1}} + 
                          \frac{s_{l-j-1}}{d_{l-1}2^{j+1}}} &
                                j < l \leq n-i+1   \end{array} \right. \]
More concretely, $a_{2,l}$ is $+1$ whenever $s_{l-j} =
s_{l-j-1}$, and is $-1$ otherwise.

By construction, 
$\theta_{1,n-i+1}({\bf a_1})$ and $\theta_{1,n-i+1}({\bf a_2})$ 
each have length no greater than $4(n-i+1)$.  
The vectors were also constructed so that 
\[\theta_{1,l}(a_{1,1}, \ldots, a_{1,l}) \theta_{1,l}(a_{2,1}, \ldots, 
a_{2,l}) = 1\]
for each $l \leq j$.  When $l > j$, we use
Proposition~\ref{thetaprops} to compute 
\[ \begin{array}{l}
\theta_{1,l}(a_{1,1}, \ldots, a_{1,l}) \theta_{1,l}(a_{2,1}, \ldots, 
a_{2,l}) \\
\hspace{.25in}=
  \theta_{1,l-1}(a_{1,1}, \ldots, a_{1,l-1}) \theta_{1,l-1}(a_{2,1}, 
\ldots, a_{2,l-1})
    \epsb{l}\left(d_l \left(2^j - 2^j(-1)^{\frac{s_{l-j}}{d_l2^{j+1}}} 
  \right) \right) \\
\hspace{.25in}=
  \theta_{1,l-1}(a_{1,1}, \ldots, a_{1,l-1}) \theta_{1,l-1}(a_{2,1}, 
\ldots, a_{2,l-1})
    \epsb{l}(s_{l-j}).\end{array}\]
By induction 
$\theta_{1,n-i+1}({\bf a_1}) 
\theta_{1,n-i+1}({\bf a_2})=\prod_{l=j+1}^{n-i+1} \epsb{l}(s_{l-j})$, 
which proves the claim.
\end{proof}

We now consider the case where ${\bf u_j}$ is zero in all but the last
$r_i-n+i-1$ coordinates.  This can only happen in types $B_n$, $C_n$, and
$D_n$.

\begin{lemma} 
Fix $j \leq r_i -1$ 
and denote $\max{\{j,n-i+1}\}$ by $q$.
Let ${\bf u} \in K^{r_i-q}$ be a vector whose entries have
each $s_l \in \{0,d_{l+q}2^{j+1}\}$.
The word $\prod_{l=q+1}^{r_i} \epsb{l}(s_{l-q})$
has length at most $24r_i$.
\label{aftermiddle}
\end{lemma}

\begin{proof}
Define vectors ${\bf a_1}$ and ${\bf a_2}$ by
\[a_{1,l} = \left\{ \begin{array}{ll}
                        2 & 1 \leq l \leq j-n+i-1 \\
                        1 &  j-n+i-1 < l \leq r_i-n+i-1 
                \end{array} \right. \]
(in particular,  if $j \leq n - i + 1$ then $a_{1,l} = 1$ for all $l$) and
\[a_{2,l} = \left\{ \begin{array}{ll}
                         2 &  1 \leq l \leq j-n+i-1 \\
                        \left( -1 \right)^{\frac{s_{l-q}}{d_{l} 2^{j+1}} + 
                          \frac{s_{l-q-1}}{d_{l-1}2^{j+1}}} &
                                j-n+i-1 < l \leq r_i-n+1-1.   
                 \end{array} \right. \]

Write $t=d_{n-i+1} 2^j$ if $j< n-i+1$ and 
$t = d_{n-i+1} 2^{n-i+1}$ otherwise.
The previous lemma showed that $\epsb{n-i+1}(\pm t)$ 
is a word of length at most $8(n-i+1)$.

By construction
$\theta_{n-i+1,r_i-n+i}(t,{\bf a_1})$ and 
$\theta_{n-i+1,r_i-n+i}(-t,{\bf a_2})$ 
each have length at most $2\left( 8(n-i+1)+4(r_i-n+i-1) \right)$.  Since
$r_i \geq 2(n-i)$ in types $B_n$, $C_n$, and $D_n$, this length 
is at most $24r_i$.  
As in the previous proof,
\[\theta_{n-i+1,r_i-n+i}(t,{\bf a_1}) \theta_{n-i+1,r_i-n+i}(-t,{\bf a_2}) =
   \prod_{l=q+1}^{r_i} \epsb{l}(s_{l-q}).\]
\end{proof}
 
The words $\prod_{l=k+1}^{r_i} \epsb{l}(s_{j,l})$ can be constructed
by Lemmas \ref{beforemiddle} and \ref{aftermiddle} as long as $k \geq
j$.  Since $j$ ranges only up to $\lfloor \frac{\log p}{\log 2}
\rfloor$, we have shown that if we take $k = \lfloor \frac{\log p}{\log 2}
\rfloor$, each word $\prod_{l=k+1}^{r_i} \epsb{l}(s_{j,l})$ can be
expressed by a word of length at most $32 r_i$.  We thus obtain:

\begin{proposition}\label{boundingU>k}
Every element of $U_i^{>k}$ can be expressed by a word in $\Sigma$ of
length at most 
$32r_i \left( \frac{\log{p}}{\log{2}} + 1 \right)$.
\end{proposition}

\subsection{Generating the kernel of the projection to $U_i$}
\label{vi}
In type $A_n$ the groups $U_i$ and $U_i'$ are identical.  In the other
types, the kernel $V_i$ of the projection from $U_i'$ to $U_i$ is
a single root subgroup $\epsb{}(\mathbb{Z}/p \mathbb{Z})$.  In order
to describe completely the generation of $U_i'$ by words in $\Sigma$,
it remains only to treat this subgroup.

We now prove that  this root subgroup can be generated using the method
of Section \ref{first k} 
if $n-i+1 \leq \lfloor \frac{\log p}{\log 2} \rfloor$ and that 
of Section \ref{last k} if not.

Let $\beta$ be the root whose root subgroup is $V_i$.
In type $B_n$ the root $\beta$ is $\sum_{j=i}^{n} \alpha_j$, while in 
type $C_n$ the root $\beta$ is $\sum_{j=i}^{n-1} 2 \alpha_j + \alpha_n$, and
in type $D_n$ the root $\beta$ is $\sum_{j=i}^{n-2} \alpha_j + \alpha_n$.

We will show below that the subgroup $V_i$ can be generated by 
explicitly given commutators.  
Define an integer $d$ by
\[d = \left\{ \begin{array}{ll}
2d_{n-i}c_{11} & \textup{ in type }B_n, \\
2d_{n-i}^2 c_{21} & \textup{ in type }C_n, \textup{ and} \\
d_{n-i-1}c_{11} & \textup{ in type }D_n \\
\end{array} \right.  \]
where the structure constants $c_{ij}$ are computed with respect to the
commutator $[\epsb{n-i}(1), \epsa{n}(1)]$ in type $B_n$ and $C_n$, and
with respect to  $[\epsb{n-i-1}(1), \epsa{n}(1)]$ in type $D_n$, and
$d_l$ is the product of the first $l$ structure constants as defined
above.  Note that $d$ is a power of $2$ and in particular is a unit in $\Z/pZ$.

\begin{lemma}
Every element of $V_i$ can be written as a word in $\Sigma$ of length
at most $32(n-i+1)(\sqrt{p}+1)$
\end{lemma}

\begin{proof}
As in the proof of Lemma \ref{rootelements}, we generate an arbitrary element
$\epsb{}(sd)$ by means of commutators.  As there, set $t=\lfloor
\sqrt{s} \rfloor$.  By contrast to Lemma \ref{rootelements}, our
commutator identities here lie in $U'_i$, not the quotient $U_i$.  We
will use the fact that the root subgroups attached to $\beta$ and
$\beta_{n-i+1}$ commute with each other, as can be seen by another
application of the Chevalley commutator relation.

In type $B_n$ we have that 
\[[\epsb{n-i}(ud_{n-i}), \epsa{n}(v)] = \epsb{}\left(uv\frac{d}{2}\right)
\epsb{n-i+1}(uv^2d_{n-i+1})\] 
and 
\[[\epsb{n-i}(-ud_{n-i}), \epsa{n}(-v)] = \epsb{}\left(uv\frac{d}{2}\right)
\epsb{n-i+1}(-uv^2d_{n-i+1}).\]  
It follows that
\[\begin{array}{rcl}
\epsb{}(t^2d)&=&[\epsb{n-i}(td_{n-i}), \epsa{n}(t)][\epsb{n-i}(-td_{n-i}), 
\epsa{n}(-t)] \textup{ and} \\
\epsb{}((s-t^2)d) &=& [\epsb{n-i}(d_{n-i}), 
\epsa{n}(s-t^2)][\epsb{n-i}(-d_{n-i}), \epsa{n}(t^2-s)].\end{array}\]
Similarly, in type $C_n$ we have that
\[[\epsb{n-i}(ud_{n-i}), \epsa{n}(v)] = \epsb{}\left(u^2v\frac{d}{2}\right)
\epsb{n-i+1}(uvd_{n-i+1})\] 
and so the following relations hold:
\[\begin{array}{rcl}
\epsb{}(t^2d) &=& [\epsb{n-i}(td_{n-i}), \epsa{n}(1)] [\epsb{n-i}(-td_{n-i}),
    \epsa{n}(1)]  \textup{ and}\\
\epsb{}((s-t^2)d) &=& 
   [\epsb{n-i}(d_{n-i}), \epsa{n}(s-t^2)] [\epsb{n-i}(-d_{n-i}),
    \epsa{n}(s-t^2)].
\end{array}\]
In type $D_n$ we use the relations:
\[\begin{array}{rcl}
\epsb{}(t^2d) &=& [\epsb{n-i-1}(td_{n-i-1}), \epsa{n}(t)]  \textup{ and}\\
\epsb{}((s-t^2)d) &=& 
   [\epsb{n-i-1}(d_{n-i-1}), \epsa{n}(s-t^2)].
\end{array}\]
By the same analysis as in Lemma \ref{rootelements}, the length of 
$\epsb{n-i}(td_{n-i})$ is at most $8(n-i)\sqrt{s}$ and the length of
$\epsb{n-i}(d_{n-i})$ is at most $8(n-i)$.  (Note that we are using
here the fact that the identity in Proposition~\ref{commutator
  generation} is correct in $U_i'$, not merely $U_i$, whenever $j \leq
n-i$.)  By definition, the length of $\epsa{n}(s-t^2)$
is at most $2 \sqrt{s}$ and the length of $\epsa{n}(t)$ is at most
$\sqrt{s}$.  The total length of $\epsb{}(sd)$ is at most 
$32(n-i+1)(\sqrt{s}+1)$, and since $s$ is no greater than $p$, the
lemma is proven.
\end{proof}

\begin{lemma}
If $n-i+1 > \left\lfloor \frac{\log p}{\log 2} \right\rfloor$ then
every element of $V_i$ can be written as a word in $\Sigma$ of length
at most $8(n-i+1) \left(\left\lfloor 
\frac{\log p}{\log 2} \right\rfloor + 1 \right)$.
\end{lemma}

\begin{proof}
The following formulas hold by induction on the Chevalley commutator
relation: in type $B_n$,
\[\begin{array}{l}
\theta_{1,n-i+1}(t_1, t_2,  \ldots, t_{n-i}, t_{n-i+1})
\theta_{1,n-i+1}(-t_1, t_2,  \ldots, t_{n-i}, -t_{n-i+1}) \\
\hspace{1in}  = \epsb{}(d \prod_{l=1}^{n-i+1}t_l), \end{array}\]
 in type $C_n$,
\[\begin{array}{l}
\theta_{1,n-i+1}(t_1, t_2,  \ldots, t_{n-i+1})
\theta_{1,n-i+1}(-t_1, t_2,  \ldots, t_{n-i+1}) \\
\hspace{1in}  = \epsb{}(dt_{n-i+1} \prod_{l=1}^{n-i}t_l^2), \end{array}\]
and in type $D_n$,
\[\begin{array}{l}
\theta_{1,n-i+1}(t_1, t_2, \ldots, t_{n-i-1},0,t_{n-i+1})
\theta_{1,n-i+1}(-t_1, t_2, \ldots, t_{n-i-1},0,0) \\
\hspace{1in}   = \epsb{}(dt_{n-i+1} \prod_{l=1}^{n-i-1}t_l).\end{array}\]
To generate $\epsb{}(2^jd)$ in types $B_n$ and $D_n$, let $t_l=2$ for
$l=1$ through $j$ and $t_l=1$ otherwise.  To generate $\epsb{}(2^{2j+1}d)$
in type $C_n$ let $t_l=2$ for $l=1$ through $j$, let $t_{n-i+1}=2$, and
let $t_l=1$ otherwise.  To generate $\epsb{}(2^{2j}d)$ in type $C_n$,
let $t_l=2$ for $l=1$ through $j$, and $t_l=1$ otherwise.

This shows that $\epsb{}(2^jd)$ can be written as a word of length at
most $8(n-i+1)$ in types $B_n$, $C_n$, and $D_n$ for each $j$ between
$0$ and $\lfloor \frac{\log p}{\log 2} \rfloor$.  Thus, by means of
binary expansion, we can write $\epsb{}(sd)$ for any $s \in \Z/p\Z$ as
a word of length at most $8(n-i+1) \left(\lfloor \frac{\log p}{\log 2} \rfloor + 1\right)$.
\end{proof}

\subsection{The diameter of $U$}

\begin{theorem}
There exists a constant $C$, independent of $n$ and $p$,  such that the
diameter of $U$ in $\Sigma$ is at most
\beq
C(np + n^2 \log p).
\eeq

\end{theorem}

\begin{proof}
We have already observed that $U$ is the product of the $U'_i$, so it
suffices to show that every element $u$ of $U_i'$ can be written as a
short word in $\Sigma$.

The results of Sections~\ref{abel}, \ref{first k} and \ref{last k} show that we
can find an element $\tilde{u}$ of $U_i'$ such that $\tilde{u}$ and
$u$ project to the same element of $U_i$, and $\tilde{u}$ is a word in
$\Sigma$ of length at most
\beq
\frac{p}{2} + \frac{48}{\log^2 2}\sqrt{p}\log^2 p + 32 r_i \left (\frac{\log p}{\log
  2} + 1 \right).
\eeq
Since $\tilde{u}u^{-1}$ lies in $V_i$, we have from Section~\ref{vi}
that $\tilde{u}u^{-1}$ can be expressed by a word of length
$32(n-i+1)(\sqrt{p}+1)$ in general, and of length $8(n-i+1)\left(\frac{\log
  p}{\log 2}+1\right)$ once $n-i+1 > \left\lfloor \frac{\log p}{\log 2}
\right\rfloor$.  In either case, $u$ can be written as a word of
length at most
\beq
\frac{p}{2} + \frac{48}{\log^2 2}\sqrt{p}\log^2 p + \frac{32}{\log
  2}\log p (\sqrt{p} + 1) +  32 r_i\left(\frac{\log p}{\log
  2} + 1 \right) + 8(n-i+1)\left(\frac{\log p}{\log 2} + 1\right).
\eeq
Now the diameter of $U$ in $\Sigma$ is bounded by the sum of the above
expression as $i$ ranges from $1$ to $n$; this sum is evidently
bounded by
\beq
\frac{np}{2} + n\left[\frac{48}{\log^2 2}\sqrt{p}\log^2 p + \frac{32}{\log
  2}\log p (\sqrt{p} + 1)\right] + {36n^2 + 4n}\left(\frac{\log p}{\log 2}
+1 \right)   
\eeq
using the fact that $\sum_{i=1}^n r_i$ is at most $n^2$ by
\cite[section 12.2]{H2}.  This expression is evidently bounded above
by a constant multiple of $np + n^2 \log p$.
 \end{proof}

\section{A lower bound for the diameter}

In this section we prove the following lower bound on the diameter of
$G$ in $\Sigma$. 

\begin{proposition} There exists an absolute constant $c$ such that, if $p$
  is an odd prime, and $G$ is a Chevalley group of classical type over
  $Z/pZ$, the diameter of $U$ in $\Sigma$ is at least $c(np + n^2 \log
  p)$.
\label{lower bound}
\end{proposition}

The main idea is to exploit the fact that each generator in $\Sigma$
commutes with almost all the others: this means that there are many
identities among the words in $\Sigma$ of some fixed length $\ell$;
from an upper bound on the number of distinct words of length $\ell$
we can derive a lower bound on the diameter of $U$. 


Since the size of our generating set for $U$ has cardinality $2n$, there
are at most $(2n)^m$ words in $G$ of length $m$.  Since $\log |U|$ is on
order of $n^2 \log p$, this implies immediately that
\beq
d(U,\Sigma) \geq \log_n |U| = O(n^2 (\log n)^{-1} \log p).
\eeq
We will show by means of commutativity between the elements of
$\Sigma$ that the $(\log n)^{-1}$ factor can be removed.  Note also
that the diameter of $U$ is at least as great as the diameter of the
abelianization of $U$, which is on order $np$.  Together, these facts
prove Proposition~\ref{lower bound}.


\medskip

For any group $G$ with a generating set $S$, let $v_{G,S}(m)$ be the
number of elements of $G$ expressible as a word of length exactly $m$
in the elements of $S$.  We organize these values into a generating
function
\beq
V_{G,S}(t) = \sum_{m=0}^\infty v_{G,S}(m) t^m.
\eeq

We take $\tilde{U}$ to be a group defined by generators and relations
as follows: let $\tilde{U}$ be generated by a set of elements  $\tilde{\Sigma} = \set{\tilde{e}_1,
\ldots, \tilde{e}_n,\tilde{e}'_1,
\ldots, \tilde{e}'_n}$, subject to the relations
\begin{equation}
[\tilde{e}_i, \tilde{e}_j] = 
[\tilde{e}'_i, \tilde{e}_j] = 
[\tilde{e}_i, \tilde{e}'_j] =
[\tilde{e}'_i, \tilde{e}'_j]
 = 1
\label{eq:tildeurel}
\end{equation}
for each pair $i,j$ such that the root subgroups of $U$ attached to
$\alpha_i$ and $\alpha_j$ commute with each other.There is a natural
surjection $\phi: \tilde{U} \ra U$ sending $\tilde{e}_i$ to $\epsa{i}(1)$ and $\tilde{e}'_i$ to $\epsa{i}(-1)$.
It follows that
\beq
v_{U,\Sigma}(m) \leq v_{\tilde{U},\tilde{\Sigma}}(m)
\eeq
for all nonnegative integers $m$.

\begin{proposition}
\beq
v_{\tilde{U},\tilde{\Sigma}}(m) \leq 2^{n+2}(4+2\sqrt{3})^m. 
\eeq
\end{proposition}

\begin{proof}
The main tool is a result of Cartier and Foata~\cite{cart:cf}, which
shows that for any group $G$ with generating set $S = \set{s_1,
\ldots, s_r}$ whose only
relations are commutations between the $s_i$,
\beq
[V_{G,S}(t)]^{-1} = \sum_T (-t)^{|T|}
\eeq
where $T$ ranges over all subsets of $S$ whose elements
commute pairwise.  For instance, if $G$ is the free abelian group on
$S$, then $T$ ranges over all subsets of $S$, and 
\beq
V_{G,S}(t)^{-1} = (1-t)^r.
\eeq
More precisely, we may define $V_{G,S}(x_1, \ldots,
x_r)$ to be the generating function whose $x_1^{m_1} \ldots
x_r^{m_r}$ coefficient is the number of distinct words in $G$ composed
of $m_1$ copies of $s_1$, $m_2$ copies of $s_2$, and so on.  Then
\beq
[V_{G,S}(x_1, \ldots, x_r)]^{-1} = \sum_T (-1)^{|T|}\prod_{i \in T}
x_i.
\eeq
 
Let $G$ be the subgroup of $\tilde{U}$ generated by $S =
\set{\tilde{e}_1, \ldots, \tilde{e}_n}$.  We begin by bounding the
volume growth of $G$ with respect to $S$.

If $\DD$ is a Dynkin diagram, write $U(\DD)$ for the unipotent group
attached to $\DD$, write $\Sigma(\DD)$ for its standard set of generators, and write
\beq
P(\DD) =  [V_{G,S}(t)]^{-1}
\eeq
where $G \subset \tilde{U}$ and $S$ are described as above.

Recall that two generators $\epsa{i}(1)$ and $\epsa{j}(1)$ in
$\Sigma(\DD)$ commute if and only if the corresponding vertices of
$\DD$ are not connected by an edge.

The pairwise-commuting subsets of
$S$ naturally split into those which contain the $n^{th}$
generator and those which do not; this yields a recursion relation
\begin{equation}
P(A_n) = P(A_{n-1}) - t P(A_{n-2})
\label{eq:recursion}
\end{equation}
for $n \geq 3$, which implies that
\beq
P(A_n) = c_1 \left( \frac{1 + \sqrt{1-4t}}{2} \right)^n + c_2 \left( \frac{1 -
\sqrt{1-4t}}{2} \right)^n
\eeq
for some constants $c_1, c_2$ depending on $t$.  Write $u =
\sqrt{1-4t}$.

The initial terms $P(A_1)=1-t = (1/4)(3 + u^2)$ and $P(A_2)=1-2t =
(1/2)(1+u^2)$ show that
\beq
c_1 = \frac{(1+u)^2}{u},
c_2 = \frac{(1-u)^2}{u}.
\eeq
So
\begin{equation}
P(A_n) = \frac{1}{u}\left[ \left(\frac{1+u}{2}\right)^{n+2} - 
\left(\frac{1-u}{2}\right)^{n+2} \right].
\label{eq:anform}
\end{equation}
In particular, when $t < 1/4$, we have $u > 0$; it follows from
\eqref{eq:anform} that $P(A_n)(t) > 0$.  So all roots of $P(A_n)$ are
at least $1/4$.  Note that if we write $P(A_0) = 1, P(A_{-1}) = 1,
P(A_{-2}) = 0$, the recursion~\eqref{eq:recursion} is still satisfied.

Since $B_n,C_n$ have the same edges as $A_n$, we have $P(A_n) = P(B_n)
= P(C_n)$.  The polynomial $P(D_n)$ satisfies the recurrence
\eqref{eq:recursion}, but has initial values
\beq
P(D_3) = 1- 3t + t^2 = P(A_3), P(D_4) = 1 - 4t + 3t^2 - t^3 = P(A_4) -
t^3.
\eeq
It follows that
\beq
P(D_n) = P(A_n) - t^3 P(A_{n-5})
\eeq
for all $n \geq 3$, since it is true for $n=3,4$ and since $P(D_n)$
satisfies \eqref{eq:recursion}.  In terms of $u$, we get
\beq
P(D_n) = \frac{1}{8u}[f(u) + f(-u)]
\eeq
where
\beq
f(u) = (1+7u-u^2-u^3)\left(\frac{1+u}{2}\right)^n.
\eeq
Now $|1 + 7u - u^2 - u^3| \geq |1 - 7u - u^2 + u^3|$ for $0 \leq u
\leq 1$; it follows that $f(u) + f(-u) > 0$ for $0 \leq u \leq 1$.
Furthermore, if $u > 1$, then $t < 0$, and it is obvious from the
recursion~\eqref{eq:recursion} that $P(D_n)(t) > 0$ for all $n$.  We
conclude that $\sqrt{1 - 4t} \leq 0$ for all $t$ such that $P(D_n)(t)
= 0$.


We have now shown that the roots of $P(\DD)$ are all at least $1/4$.

Now define
\beq
Q(\DD)(t) =  [V_{\tilde{U},\tilde{\Sigma}}(t)]^{-1}
\eeq
and
\beq
Q(\DD)(x_1, \ldots, x_n, x'_1, \ldots x'_n) =  [V_{\tilde{U},\tilde{\Sigma}}(x_1, \ldots, x_n,x'_1,
\ldots, x'_n)]^{-1}.
\eeq
The relations~\eqref{eq:tildeurel} tell us that
\beq
Q(\DD)(x_1, \ldots, x_n, x'_1, \ldots x'_n) = P(\DD)(x_1 + x'_1 -
x_1x'_1, \ldots, x_n + x'_n -
x_nx'_n)
\eeq
from which we have 
\beq
Q(\DD)(t) = P(\DD)(2t - t^2).
\eeq
So the roots of $Q(\DD)$ are all greater than $1-(\sqrt{3}/2)$.  It
follows that 
\beq
v_{\tilde{U},\tilde{\Sigma}}(m) \leq C (4 + 2\sqrt{3})^m
\eeq
for some constant $C$.  To bound this constant roughly, write
\beq
\begin{array}{l}
C = \max_{m} (4+ 2\sqrt{3})^{-m} v_{\tilde{U},\tilde{\Sigma}}(m) \\
\hspace{.5in} \leq \sum_{m}(4+ 2\sqrt{3})^{-m} v_{\tilde{U},\tilde{\Sigma}}(m)
= [Q(\DD)(1-(\sqrt{3}/2))]^{-1} = [P(\DD)(1/4)]^{-1}
\end{array}
\eeq
We can rewrite \eqref{eq:anform} as 
\beq
\begin{array}{l}
P(A_n) = \frac{1}{2^{n+2}} \left[ \frac{1}{u}(1+u)^{n+1} + (1+u)^{n+1}
   - \frac{1}{u}(1-u)^{n+1} + (1-u)^{n+1} \right] \\
\hspace{.5in}  = \cdots = \frac{1}{2^{n+2}} \left( \sum_{i=0}^{n+1} (1+u)^i + 
   (1-u)^i \right).
\end{array}
\eeq
Setting $u=0$ in this expression (that is, evaluating $P(A_n)(t)$ at $t=1/4$)
we see that
\beq
P(A_n)(1/4) = 2^{-n-1} (n+2)
\eeq
and
\beq
P(D_n)(1/4) = P(A_n)(1/4) - (1/4)^3 P(A_{n-5})(1/4) = (n+7)2^{-n-2}.
\eeq
In either case, we have $[P(\DD)(1/4)]^{-1} \leq 2^{n+2}$, which
yields the desired result.
\end{proof}
We have now shown that
\beq
v_{U,\Sigma}(m) \leq v_{\tilde{U},\tilde{\Sigma}}(m) \leq
2^{n+2}(4+2\sqrt{3})^m.
\eeq
So
\beq
\sum_{m=0}^D v_{U,\Sigma}(m) \leq (1/3)2^{n+3}\sqrt{3}(4+2\sqrt{3})^D.
\eeq

On the other hand, if $D = d(U,\Sigma)$, it must be the case that
\beq
\sum_{m=0}^D v_{U,\Sigma}(m) = |U| \geq p^{(1/2)n^2}.
\eeq
We conclude that
\beq
d(U,\Sigma) \geq [(1/2)n^2\log(p) - (n+3)\log 2 + (1/2)\log 3] / \log(4 + 2\sqrt{3})
\eeq

which is evidently bounded below by a constant multiple of $n^2 \log
p$.

\begin{rem}  The discussion here applies equally well to any group
endowed with a set of generators obeying the same
commutation relations as $\Sigma$.  For instance, the number of words
of length $m$ in the standard $n-1$ generators of the
$n$-strand Artin braid group is of order at most $4^m$.  This suggests that
the study of random walks on the Artin braid group might resemble that
of random walks on the free group on $4$ letters; this theme is
explored in detail in~\cite{vers:braidgrowth}.
\end{rem}

\section{Appendix}

In this appendix, we give explicit matrix descriptions of the
different Chevalley groups discussed in the paper.  To begin,
in type $A_n$, the group $G$ is $SL_{n+1}(K)$; in type $B_n$, 
the group $G$ is 
$SO_{2n+1}(K)$; in type $C_n$, the group $G$ is $Sp_{2n}(K)$; and in type 
$D_n$, the group
$G$ is $SO_{2n}(K)$. In all these cases, we will describe a natural
embedding of $G$ into a general linear group $GL_N(K)$ and describe
matrix representatives for generators of the $U_\alpha$.

The group $G$ is the subgroup of $GL_N(K)$ that preserves an
alternating bilinear form in type $C_n$ and an
symmetric bilinear form in types $B_n$ and $D_n$,  so $N$ is either
$2n$ or $2n+1$. 
We use the bilinear form
defined by
$\langle e_i, e_{N+1-i} \rangle = 1$ when $1 \leq i \leq N/2$ and $\pm 1$ 
otherwise, depending on whether the form is symmetric or alternating,
to write $G$ explicitly.  We denote by $E_{ij}$ the matrix having
$ij$th entry $1$ and all other entries $0$.


The $U_{\alpha}$ may now be described in terms of these explicit matrix 
representatives.   
In type $A_n$, the subgroups
$U_{\alpha}$ correspond to the groups 
$\{\Id + cE_{ij}: c \in K, i < j\}$.  By contrast, in 
type $C$ the $U_{\alpha}$ are either of the form
\[\{\Id + c(E_{ij}-E_{2n+1-j,2n+1-i}): i<j \leq n, c \in K\}\]
or of the form
\[\{\Id + c(E_{i,2n+1-j}+E_{j,2n+1-i}): i \leq j \leq n, c \in K\}.\]
In type $B_n$ the $U_{\alpha}$ are of the form
\[\{\Id + c(E_{ij}-E_{2n+2-j,2n+2-i}): i<j<2n+2-i, j \neq n+1, i< n+1, c \in 
K\}\]
or of the form
\[\{\Id + c(E_{i,n+1}-E_{n+1,2n+2-i}) - \frac{c^2}{2}E_{i,2n+2-i}: i< n+1, 
c \in K\}.\]
The $U_{\alpha}$ in type $D_n$ are always of the form
\[\{\Id + c(E_{ij}-E_{2n+1-j,2n+1-i}): i<j<2n+1-i, i< n, c \in K\}.\]


A different bilinear form could have been used and other
generators could have been chosen.  This would only 
change the specific constants
in Figure \ref{constants}.
We chose as we did so that $U$
is a subgroup of upper triangular matrices for all types.

In type $A_n$ the group $U_i' = \{(\Id + \sum_{j=i+1}^n m_{ij} E_{ij}): m_{ij} 
\in K\}$ is the group of unipotent upper-triangular matrices with nonzero
off-diagonal entries only on the $i^{th}$ row.  In the other types, $U_i'$
also has entries on the $i^{th}$ row, but the Chevalley type of the group
requires additional nonzero entries.
For instance, in type $B_n$ the group $U_i'$ is
\[\left\{\begin{array}{l}
\Id + \sum_{j=i+1}^{2n+1-i} m_{ij} \left(E_{ij}-E_{2n+2-j,2n+2-i} \right) 
    \\ \hspace{.5in} + \left( \frac{-m_{i,n+1}^2}{2} - \sum_{j=i+1}^{n} m_{ij}m_{i,2n+2-j} 
     \right) E_{i,2n+2-i}: m_{ij} \in K \end{array} \right\}.\]
In type $C_n$ the group $U_i'$ is
\[\begin{array}{l}
\displaystyle \left\{\Id + \sum_{j=i+1}^{n} m_{ij} \left(E_{ij}-E_{2n+1-j,2n+1-i} \right)
   + \sum_{j=n+1}^{2n-i} m_{ij} \left(E_{ij}+E_{2n+1-j,2n+1-i} \right)
      \right.
     \\ \displaystyle \left. \hspace{.5in} 
     + \left( m_{i,2n+1-i} + \sum_{j=i+1}^{n} m_{ij}m_{i,2n+1-j} 
     \right) E_{i,2n+1-i}: m_{ij} \in K \right\} \end{array}\]
and in type $D_n$ the group $U_i'$ is
\[\left\{\Id + 
    \sum_{j=i+1}^{2n-i} m_{ij} \left(E_{ij}-E_{2n+1-j,2n+1-i} \right):
      m_{ij} \in K \right\}.\]
If $i$ is fixed then as $j>i$ varies, the 
root subgroups $U_{\alpha}$ given above are precisely 
the root subgroups for each $\alpha$ in 
$\Ph{i}$.  With one exception, we take $\epsa{}(c)$ for each $c$ 
to be the matrix given in our description of $U_{\alpha}$.  
For instance, if $\alpha$ is the root corresponding to a pair $i$, $j \leq n$, 
then in type $A_n$ the matrix $\epsa{}(1) = \Id + E_{ij}$ while in type $B_n$
the matrix $\epsa{}(1) = \Id + E_{ij} - E_{2n+2-j,2n+2-i}$.  
(We do not consider the root subgroups $\epsa{}(1) = \Id + E_{i,2n+1-i} + 
E_{i,2n+1-i} = \Id + 2E_{i,2n+1-i}$ in type $C_n$ to be an exception
to this rule.)
The one exception
to this rule is for the simple root $\alpha_n$ in type $D_n$.  In this 
case, we take
\[\epsa{n}(1) = \Id - E_{n,n+2}+E_{n-1,n+1},\]
and in general in type $D_n$ the matrix $\epsa{n}(c) = \Id +c(-E_{n,n+2}
+E_{n-1,n+1})$.


Fix $\alpha$ in $\Ph{i}$ and take $j \geq i$.  Assume that $\alpha + \alpha_j
\in \Phi^+$.
The commutator $[\epsa{}(s), \epsa{j}(t)]$
always lies in exactly one root subgroup in types $A_n$ and $D_n$, and
almost always in types $B_n$ and $C_n$.
Using the matrices which generate the $U_{\alpha}$ we compute the 
structure constants $c_{ij}$ explicitly.  
The following table gives the values of $c_{ij}$ in the Chevalley
commutator relation for all of the
cases we will consider.   
The row with the ordered pair of roots $(\alpha, \alpha_j)$ 
contains the nonzero coefficients in the Chevalley 
commutator expansion of $[\epsa{}(s), \epsa{j}(t)]$.  

\begin{figure}[h]
 \begin{tabular}{|c|c|c|} 

        \cline{1-3}
        Ordered pair & Nonzero coefficients & Conditions \\
        \cline{1-3} 
        $(\alpha,\alpha_j)$ & $c_{11}=1$ & $\begin{array}{l}
                           \alpha \not> \alpha_j, \\ 
                           j \neq n \textup{ in types $B_n$, $C_n$, $D_n$},\\
                                  j \neq n-1 \textup{ in type $D_n$}
                                  \end{array}$\\
        \cline{1-3} 
        $(\alpha,\alpha_j)$ & $c_{11}=-1$  & $\alpha > \alpha_j$ \\
        \cline{1-3}
        $(\sum_{j=i}^{n-1} \alpha_j,\alpha_n)$ & $c_{11} = -1, 
              c_{12}=\frac{1}{2}$
                       & $\textup{ type } B_n$ \\
        \cline{1-3}
        $(\sum_{j=i}^{n-1} \alpha_j,\alpha_n)$ & $c_{11} = -2,c_{21}=2$
                       & $\textup{ type } C_n$ \\
        \cline{1-3}
        $(\sum_{j=i}^{n-2} \alpha_j,\alpha_{n-1})$ & $c_{11} =1$
                       & $\textup{ type } D_n$ \\ 
        \cline{1-3}
        $(\sum_{j=i}^{n-2} \alpha_j + \alpha_n,\alpha_{n-1})$ & $c_{11} =-1$
                       & $\textup{ type } D_n$ \\
        \cline{1-3}
        $(\sum_{j=i}^{n-2} \alpha_j,\alpha_{n})$ & $c_{11} =1$
                       & $\textup{ type } D_n$ \\
        \cline{1-3}
        $(\sum_{j=i}^{n-1} \alpha_j,\alpha_{n})$ & $c_{11} =-1$
                       & $\textup{ type } D_n$ \\
         \hline 
          \end{tabular} 
\caption{Structure constants for $U$} \label{constants} 
\end{figure}

\end{document}